# A Family of Wavelets and a new Orthogonal Multiresolution Analysis Based on the Nyquist Criterion

H. M. de Oliveira, L. R. Soares, T. H. Falk

**Abstract** — *A generalisation of the Shannon complex wavelet is introduced, which is related to raised cosine filters. This approach is then used to derive a new family of orthogonal complex wavelets based on the Nyquist criterion for Inter-symbolic Interference (ISI) elimination. An orthogonal Multiresolution Analysis (MRA) is presented, showing that the roll-off parameter should be kept below 1/3. The pass-band behaviour of the wavelet Fourier spectrum is examined. The left and right roll-off regions are asymmetric; nevertheless the Q-constant analysis philosophy is maintained. A generalisation of the (square root) raised cosine wavelets is proposed. Finally, a computational implementation of such wavelets on Matlab$^{TM}$ is presented as well as a few applications.*

## 1. Introduction

Wavelet analysis has matured rapidly over the past years and has been proved to be valuable for both scientists and engineers [1, 2]. Wavelet transforms have lately gained prolific applications throughout an amazing number of areas [3, 4]. Another strongly related tool is the multiresolution analysis (MRA). Since its introduction in 1989 [5], MRA representation has emerged as a very attractive approach in signal processing, providing a local emphasis of features of importance to a signal [2, 6, 7]. The purpose of this paper is twofold: first, to introduce a new family of wavelets and then to provide a new and complete orthogonal multiresolution analysis. Within the broad spectrum of wavelet application in Telecommunications, two among them have special motivation: (*i*) *Wavelets* in Multiplex and coding division multiple access CDMA [8 - 11], (*ii*) *Wavelets* in orthogonal multicarrier systems (type orthogonal frequency division multiplex OFDM) [12 - 20]. It has been shown that wavelet-based OFDM system is superior regarding classical OFDM multicarrier systems in many aspects [15]. Orthogonal wavelets have also been used as spreading spectrum sequences [21, 22].

We adopt the symbol := to denote "*equals by definition*". As usual, $Sa(t):=\sin(t)/t$ is adopted in all mathematical development. Nevertheless, the MRA area is accustomed to $Sinc(t):=\sin(\pi t)/\pi t$ notation (*Sinc*-MRA). The gate function of length T is denoted by $\prod\left(\frac{t}{T}\right)$ and $\delta(t)$ denotes the Dirac distribution. Wavelets are denoted by $\psi(t)$ and scaling functions by $\phi(t)$. The paper is organised as follows. Section 2 generalises the *Sinc*-MRA. A new orthogonal MRA based on the raised cosine is introduced in Section 3. A new family of orthogonal wavelet is also given. These new orthogonal wavelets seem to be particularly suitable to replace *Sinc*(.) pulses in OFDM systems. Further generalisations are carried out in Section 4. A computational implementation for analysing signals by Matlab$^{TM}$ wavelet toolbox is described in Section 5 and an example is given. Finally, Section 6 presents conclusions and perspectives.

## 2. A Generalised Shannon Wavelet (the Raised-Cosine Wavelet)

The scaling function for the Shannon MRA (or *Sinc*-MRA) is given by the "sample function": $\phi^{(Sha)}(t) = Sinc(t)$. A naive and interesting generalisation of the complex Shannon wavelet can be done by using spectral properties of the raised-cosine filter [23]. The most common filter in Digital Communication Systems, the raised cosine spectrum $P(w)$ with a roll-off factor $\alpha$, was conceived to eliminate the **Inter-symbolic Interference** (ISI). Its transfer function is given by

$$P(w) = \begin{cases} \frac{1}{2\pi} & 0 \leq |w| < (1-\alpha)\pi \\ \frac{1}{4\pi}\left\{1 + \cos\frac{1}{2\alpha}(|w| - \pi(1-\alpha))\right\} & (1-\alpha)\pi \leq |w| < (1+\alpha)\pi \\ 0 & |w| \geq (1+\alpha)\pi. \end{cases} \quad (1)$$

The "raised cosine" frequency characteristic therefore consists of a flat spectrum portion followed by a roll-off portion with a sinusoidal format. Such spectral shape is very often used in the design of base-band digital systems. It is derived from the pulse shaping design criterion that would yield zero ISI, the so-called Nyquist Criterion. Note that $P(w)$ is a real and non-negative function [23], and in addition

$$\sum_{l \in Z} P(w + l2\pi) = \frac{1}{2\pi}. \quad (2)$$

Furthermore, the following normalisation condition holds: $\frac{1}{2\pi}\int_{-\infty}^{+\infty} P(w)dw = 1$. It is proposed at this point the replacement of the Shannon scaling function on the frequency domain by a raised cosine, with parameter $\alpha$ (Fig. 1). Assume that $\Phi(w)=P(w)$. In the time domain this corresponds exactly to the impulse response of a Nyquist raised-cosine filter.

Authors are with Federal University of Pernambuco, UFPE, C.P.7800, 50711-970, Recife-PE, Brazil. E-mail: {hmo,lusoares}@ufpe.br, t_falk@terra.com.br. This paper is dedicated to Dr. Max Gerken (*in memoriam*).

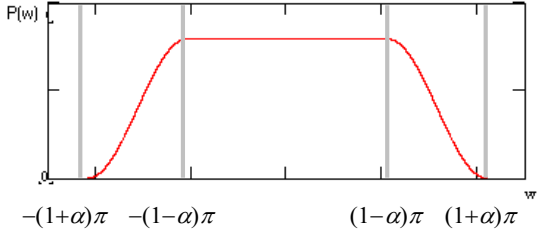

$-(1+\alpha)\pi \quad -(1-\alpha)\pi \qquad (1-\alpha)\pi \quad (1+\alpha)\pi$

**Figure 1**. Fourier Spectrum of the raised cosine scaling function: a flat spectrum portion followed by two roll-off symmetrical portions with a sinusoidal format.

The generalised Shannon scaling function is therefore:

$$\phi^{(GSha)}(t) := \frac{\cos\alpha\pi t}{1-(2\alpha t)^2} Sinc(t) \cdot \quad (3)$$

In the particular case $\alpha=0$, the scaling function simplifies to the classical Shannon scaling function. As a consequence of the Nyquist criterion, the scaling function presents zero crossing points on the unidimensional grid of integers, $n=\pm 1, \pm 2, \pm 3,\ldots$. This scaling function $\phi$ defines a non-orthogonal MRA. Figure 2 shows the scaling function corresponding to a generalised Shannon MRA for a few values of $\alpha$.

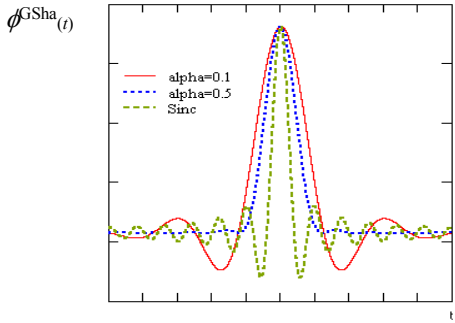

**Figure 2**. Scaling function for the raised cosine wavelet (generalised Shannon scaling function for $\alpha = 0.1$ and 0.5. The Sinc function is also plotted for comparison purposes).

## 3. Multiresolution Analysis Based on Nyquist Filters

A very simple way to build an orthogonal MRA via the raised cosine spectrum [23] can be accomplished by invoking Meyer's central condition [6]:

$$\sum_{n\in Z}|\Phi(w+2\pi n)|^2 = \frac{1}{2\pi} \cdot \quad (4)$$

Comparing Eqn(2) to Eqn(4), we choose $\Phi(w) = \sqrt{P(w)}$ (i.e. a square root of the raised cosine spectrum). Let then

$$\Phi(w)=\begin{cases} \frac{1}{\sqrt{2\pi}} & 0\le |w|<(1-\alpha)\pi \\ \frac{1}{\sqrt{2\pi}}\cos\frac{1}{4\alpha}(|w|-(1-\alpha)\pi) & (1-\alpha)\pi \le |w|<(1+\alpha)\pi \\ 0 & |w|\ge (1+\alpha)\pi. \end{cases} \quad (5)$$

Clearly, $\sum_n |\Phi(w+2\pi n)|^2 = \frac{1}{2\pi}$, so the square root of the raised cosine shape allows an orthogonal MRA. The scaling function $\phi(t)$ is plotted in the spectral domain (Fig. 3).

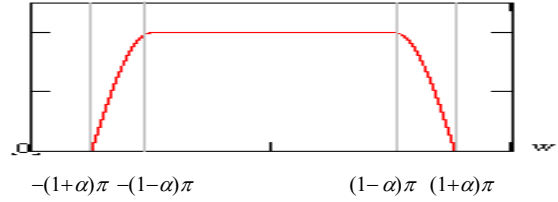

$-(1+\alpha)\pi \quad -(1-\alpha)\pi \qquad (1-\alpha)\pi \quad (1+\alpha)\pi$

**Figure 3**. Spectral Characteristic of the "de Oliveira" scaling function for an orthogonal MRA: Despite the slight shape difference as compared to Fig.1, the former can just perform an non-orthogonal MRA while this scaling function achieves an orthogonal MRA.

The cosine pulse function PCOS defined below plays an important role on the raised cosine MRA.

**Definition 1**. The cosine pulse function of parameters ($t_0$, $\theta_0$, $w_0$) and $B$ is defined by

$$PCOS(w;t_0,\theta_0,w_0,B):=\cos(wt_0+\theta_0)\prod\left(\frac{w-w_0}{2B}\right), t_0,\theta_0,w_0,B\in\Re,$$

$0<B<w_0$. □

It corresponds to a cosine pulse (in the frequency domain), with frequency $t_0$ and phase $\theta_0$, with duration $2B$ rad/s, centred at $w_0$ rad/s. Some interesting particular cases include:

1) The Gate function: $\prod\left(\frac{w}{2B}\right) = PCOS(w;0,0,0,B)$,

2) A Gate shifted by $w_0$ $\prod\left(\frac{w-w_0}{2B}\right) = PCOS(w;0,0,w_0,B)$,

3) An infinite cosine pulse:
$$\cos(wt_0+\theta_0) = PCOS(w;t_0,\theta_0,0,B\to+\infty).$$

Denoting the inverse Fourier transform of PCOS by $pcos(t;t_0,\theta_0,w_0,B):=\Im^{-1} PCOS(w;t_0,\theta_0,w_0,B)$, the following result can be proved.

**Proposition 1**. Given $t_0$, $\theta_0$, $w_0$ and $B$ parameters of a PCOS, the inverse spectrum pcos is given by:

$$pcos(t;t_0,\theta_0,w_0,B) =$$
$$\frac{B}{2\pi}\left\{e^{j(w_0t+w_0t_0+\theta_0)}.Sa[B(t+t_0)] + e^{j(w_0t-w_0t_0-\theta_0)}.Sa[B(t-t_0)]\right\} \quad (6)$$

**Proof**. Follows applying the convolution property for the following couple of transform pairs:

$$\frac{1}{2}\left[\delta(t+t_0)e^{j\theta_0} + \delta(t-t_0)e^{-j\theta_0}\right] \leftrightarrow \cos(wt_0+\theta_0) \quad , \text{ and}$$

$$\frac{B}{\pi}e^{jw_0t}.Sa(Bt) \leftrightarrow \prod\left(\frac{w-w_0}{2B}\right). \quad \square$$

It is interesting to check some particular cases:

$pcos(w;0,0,0,B) \leftrightarrow PCOS(w;0,0,0,B) \Leftrightarrow \frac{B}{\pi}.Sa(Bt) \leftrightarrow \prod\left(\frac{w}{2B}\right)$

$pcos(t;t_0,0,0,B\to+\infty) \leftrightarrow PCOS(w;t_0,0,0,B\to+\infty)$

$\Leftrightarrow \frac{1}{2}[\delta(t+t_0)+\delta(t-t_0)] \leftrightarrow \cos(wt_0)$, which follows from the property of the sequence

$$\lim_{\varepsilon\to 0} \frac{1}{\pi\varepsilon}Sa\left(\frac{t}{\varepsilon}\right) = \delta(t) \cdot \quad (7)$$

**Property 1**. (Time shift): A shift T in time is equivalent to the following change of parameters: $pcos(t-T;t_0,\theta_0,w_0,B) = pcos(t;t_0-T,\theta_0,w_0,B)$. □

## 3.1 Scaling Function Derived from Nyquist Filters

In order to find out the scaling function of the new orthogonal MRA introduced in this section, let us take the inverse Fourier transform of $\Phi(w)$.

The spectrum $\Phi(w)$ can be rewritten as a sum of contributions from three different sections (a central flat section and two cosine-shaped ends):

$$\sqrt{2\pi}\Phi(w) = \prod\left(\frac{w}{2\pi - 2B}\right) + cos(wt_0 + \theta_0)\prod\left(\frac{w - \pi}{2B}\right) + \quad (8)$$
$$cos(-wt_0 + \theta_0)\prod\left(\frac{-w - \pi}{2B}\right),$$

with parameters $B := \pi\alpha$, $t_0 := \frac{1}{4\alpha}$ and $\theta_0 := -\frac{(1-\alpha)\pi}{4\alpha}$.

It follows from Definition 1 that
$$\sqrt{2\pi}\Phi(w) = PCOS(w;0,0,0,2\pi - 2B) +$$
$$PCOS\left(w; \frac{1}{4\alpha}, -\frac{(1-\alpha)\pi}{4\alpha}, \pi, \pi\alpha\right) + PCOS\left(-w; \frac{1}{4\alpha}, -\frac{(1-\alpha)\pi}{4\alpha}, \pi, \pi\alpha\right)$$

and therefore
$$\sqrt{2\pi}\phi^{(deO)}(t) = pcos(t;0,0,0,2\pi - 2B) +$$
$$pcos\left(t; \frac{1}{4\alpha}, -\frac{(1-\alpha)\pi}{4\alpha}, \pi, \pi\alpha\right) + pcos\left(-t; \frac{1}{4\alpha}, -\frac{(1-\alpha)\pi}{4\alpha}, \pi, \pi\alpha\right).$$

After a somewhat tedious algebraic manipulation, we derive

$$\phi^{(deO)}(t) = \frac{1}{\sqrt{2\pi}} \cdot (1-\alpha) \cdot Sinc[(1-\alpha)t] + \quad (9)$$
$$\frac{1}{\sqrt{2\pi}} \cdot \frac{4\alpha}{\pi} \cdot \frac{1}{1-(4\alpha t)^2}\{cos\pi(1+\alpha)t + 4\alpha t \cdot sin\pi(1-\alpha)t\}.$$

A sketch of the above orthogonal MRA scaling function is shown in Figure 4, assuming a few roll-off values.

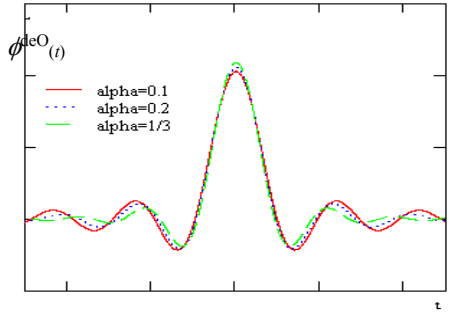

**Figure 4**. "*de Oliveira*" Scaling Function for an Orthogonal Multiresolution Analysis ($\alpha$=0.1, 0.2 and 1/3).

The scaling function $\phi^{(deO)}(t)$ can be expressed in a more elegant and compact representation with the help of the following special functions:

**Definition 2.** (Special functions); $\nu$ is a real number, $H_\nu(t) := \nu Sinc(\nu t)$, $0 \leq \nu \leq 1$, and

$$M_{\nu_2}^{\nu_1}(t) := \frac{1}{\pi}\frac{2|\nu_1 - \nu_2|}{1-[2t(\nu_1 - \nu_2)]^2}\{cos\pi\nu_1 t + 2(\nu_1 - \nu_2)t \cdot sin\pi\nu_2 t\}$$
□

It follows that:
$\sqrt{2\pi}\phi^{(Sha)}(t) = H_1(t)$,
$\sqrt{2\pi}\phi^{(deO)}(t) = H_{1-\alpha}(t) + M_{1-\alpha}^{1+\alpha}(t)$.
Clearly, $\lim_{\alpha \to 0} \phi^{(deO)}(t) = \phi^{(Sha)}(t)$.

The low-pass $H(.)$ filter of the MRA can be found by using the so-called two-scale relationship for the scaling function [7]:

$$\Phi(w) = \frac{1}{\sqrt{2}}H\left(\frac{w}{2}\right)\Phi\left(\frac{w}{2}\right). \quad (10)$$

How should $H$ be chosen to make Eqn(10) hold? Initially, let us sketch the spectrum of $\Phi(w)$ and $\Phi(w/2)$ as shown in Figure 5.

The main idea is to not allow overlapping between the roll-off portions of these spectra. Imposing that $2\pi(1-\alpha)>(1+\alpha)\pi$, it follows that $\alpha<1/3$ (remember that $0<\alpha<1$). This is a simplifying hypothesis. It is quite usual the use of small roll-off factors in Digital Communication Systems.

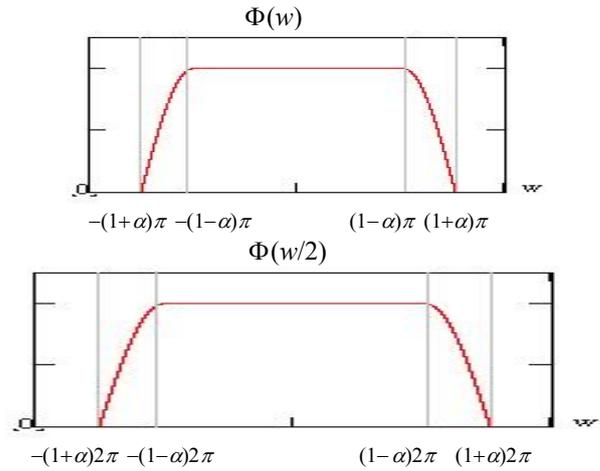

**Figure 5**. Draft of $\Phi(w)$ and $\Phi(w/2)$.

It is suggested to assume that $H\left(\frac{w}{2}\right) = \frac{1}{\sqrt{\pi}}\Phi(w)$. Substituting this transfer function into the refinement equation (Eqn(10)), results in

$$\Phi(w) = \frac{1}{\sqrt{2}}\frac{1}{\sqrt{\pi}}\Phi(w)\Phi\left(\frac{w}{2}\right). \quad (11)$$

The above equation is actually an identity for $|w| > (1+\alpha)\pi$. Into the region $|w| < (1+\alpha)\pi$, it can be seen that $\Phi\left(\frac{w}{2}\right) = \sqrt{2\pi}$, under the constraint $\alpha<1/3$.

## 3.2 The Orthogonal "*de Oliveira*" Wavelet

The orthogonal "*de Oliveira*" wavelet can be found by the following procedure [7]:

$$\Psi(w) = e^{-jw/2}\frac{1}{\sqrt{2}}H^*\left(\frac{w}{2} - \pi\right)\Phi\left(\frac{w}{2}\right). \quad (12)$$

Inserting the shape of the filter $H$ in the previous equation, it follows that:

$$\Psi(w) = e^{-jw/2}\frac{1}{\sqrt{2\pi}}\Phi(w - 2\pi)\Phi\left(\frac{w}{2}\right). \quad (13)$$

In order to evaluate the spectrum of the mother wavelet, we plot both $\Phi(w - 2\pi)$ and $\Phi\left(\frac{w}{2}\right)$, again under the constraint $\alpha<1/3$ (Fig. 6). In this case, $(1+\alpha)\pi<(1-\alpha)2\pi$ and $(1+\alpha)2\pi<(1-\alpha/3)3\pi$.

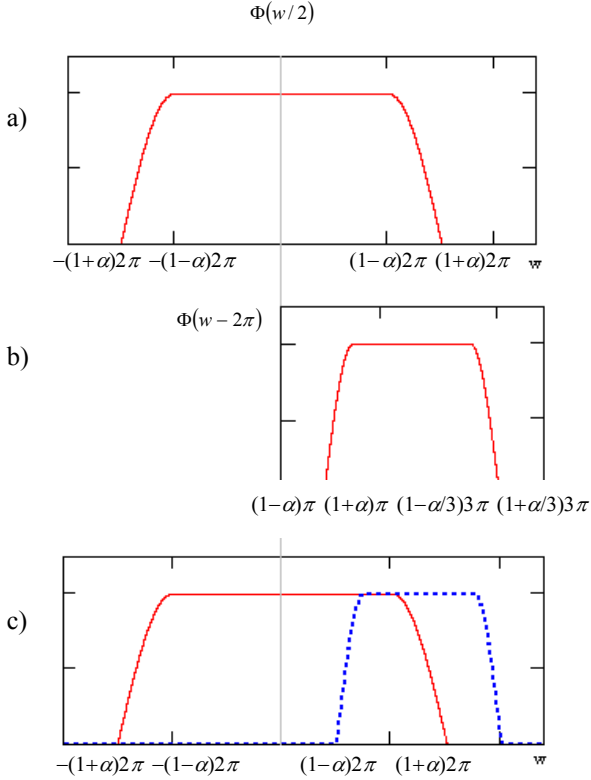

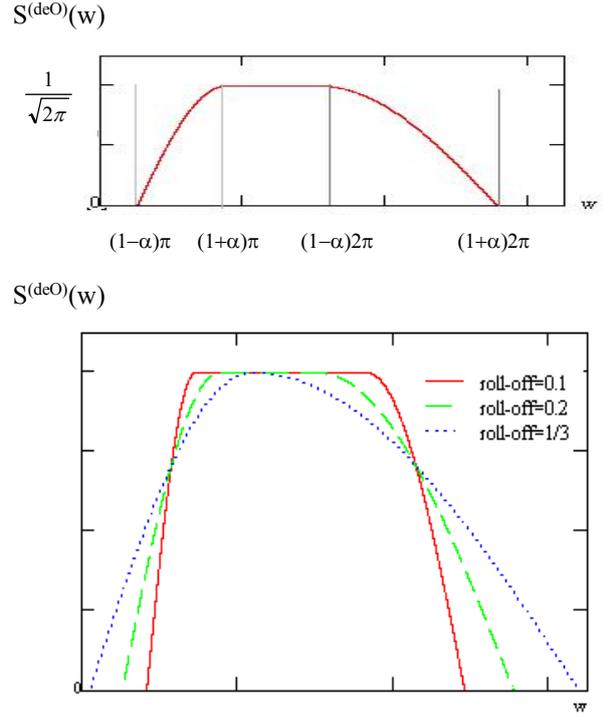

**Figure 6**. Sketch of the scaling function spectrum: (a) scaled version, $\Phi\left(\frac{w}{2}\right)$ (b) translated version $\Phi(w-2\pi)$, (c) simultaneous plot of (a) and (b).

Defining a shaping pulse

$$S^{(deO)}(w) := \frac{1}{\sqrt{2\pi}}\Phi(w-2\pi)\Phi\left(\frac{w}{2}\right), \quad (14)$$

the wavelet specified by Eqn(12) can be rewritten as $\Psi^{(deO)}(w) = e^{-jw/2}S^{(deO)}(w)$. The term $e^{-jw/2}$ accounts for the *wave*, while the term $S(w)$ accounts for *let*.

From the Figure 6, it follows by inspection that:

$$S^{(deO)}(w) = \begin{cases} 0 & \text{if } w < \pi(1-\alpha) \\ \Phi(w-2\pi) & \text{if } \pi(1-\alpha) \leq w < \pi(1+\alpha) \\ \frac{1}{\sqrt{2\pi}} & \text{if } \pi(1+\alpha) \leq w < 2\pi(1-\alpha) \\ \Phi\left(\frac{w}{2}\right) & \text{if } 2\pi(1-\alpha) \leq w < 2\pi(1+\alpha) \\ 0 & \text{if } w \geq 2\pi(1+\alpha). \end{cases} \quad (15)$$

Inserting the (square root) raised cosine format of $\Phi(.)$, results in:

$$S^{(deO)}(w) = \begin{cases} 0 & \text{if } w < \pi(1-\alpha) \\ \frac{1}{\sqrt{2\pi}}\cos\frac{1}{4\alpha}(w-\pi(1+\alpha)) & \text{if } \pi(1-\alpha) \leq w < \pi(1+\alpha) \\ \frac{1}{\sqrt{2\pi}} & \text{if } \pi(1+\alpha) \leq w < 2\pi(1-\alpha) \\ \frac{1}{\sqrt{2\pi}}\cos\frac{1}{8\alpha}(w-2\pi(1-\alpha)) & \text{if } 2\pi(1-\alpha) \leq w < 2\pi(1+\alpha) \\ 0 & \text{if } w \geq 2\pi(1+\alpha). \end{cases}$$

(16)

The complex "*de Oliveira*" wavelet is given by $\Psi^{(deO)}(w) = e^{-jw/2}S^{(deO)}(w)$ and its modulo $|\Psi^{(deO)}(w)| = S^{(deO)}(w)$ is depicted in Figure 7. Observe furthermore that making $\alpha \to 0$, the wavelet reduces to the complex Shannon wavelet.

**Figure 7**. Shaping pulse of the "*de Oliveira*" Wavelet (frequency domain, $(1-\alpha)\pi \leq |w| \leq (1+\alpha)2\pi$). The magnitude of the flat central portion is $1/\sqrt{2\pi}$.

It is quite apparent from Fig. 7 the band-pass behaviour of the wavelet $\Psi^{(deO)}(w)$. Observe that the left and right roll-off is not exactly symmetrical. Instead, despite their similar shape, they occur at different scales, a typical behaviour of wavelets.

Taking the inverse Fourier transform can derive the time domain representation of the wavelet: $\psi^{(deO)}(t) = \mathfrak{F}^{-1}\Psi^{(deO)}(w)$.

Denoting by $s^{(deO)}(t) \leftrightarrow S^{(deO)}(w)$ the corresponding transform pair, it follows that $\psi^{(deO)}(t) = s^{(deO)}\left(t-\frac{1}{2}\right)$.

The shaping pulse can be rewritten as:

$$\sqrt{2\pi}S^{(deO)}(w) = PCOS\left(w; \frac{1}{4\alpha}, -\frac{\pi(1+\alpha)}{4\alpha}, \pi, \pi\alpha\right) + PCOS\left(w; 0, 0, \frac{3\pi}{2}\left(1-\frac{\alpha}{3}\right), \frac{\pi}{2}(1-3\alpha)\right) + PCOS\left(w; \frac{1}{8\alpha}, -\frac{2\pi(1-\alpha)}{8\alpha}, 2\pi, 2\pi\alpha\right)$$

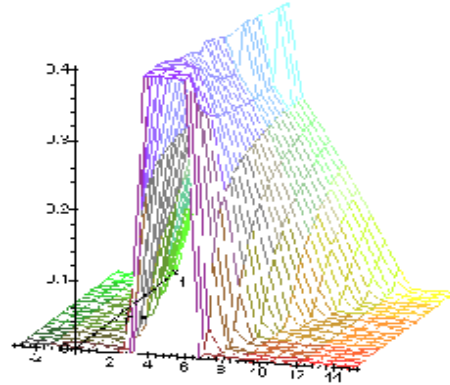

**Figure 8**. Modulo of "*de Oliveira*" Wavelet on frequency domain varying the roll-off parameter (depth axis).

Finally, applying the inverse transform, we have

$$\sqrt{2\pi}s^{(deO)}(t) = pcos\left(t;\frac{1}{4\alpha},-\frac{\pi(1+\alpha)}{4\alpha},\pi,\pi\alpha\right) + \quad (17)$$

$$pcos\left(t;0,0,\frac{3\pi}{2}(1-\frac{\alpha}{3}),\frac{\pi}{2}(1-3\alpha)\right) + pcos\left(t;\frac{1}{8\alpha},-\frac{2\pi(1-\alpha)}{8\alpha},2\pi,2\pi\alpha\right)$$

The *pcos*(.) signal is a complex signal when there are no symmetries in *PCOS*(.). The real and imaginary parts of the *pcos* function can be handled separately, according to

$pcos(t;t_0,\theta_0,w_0,B) = rpcos(t) + j.ipcos(t)$, where

$rpcos(t) := \Re(pcos(t;t_0,\theta_0,w_0,B))$ and

$ipcos(t) := \Im(pcos(t;t_0,\theta_0,w_0,B))$.

Aiming to investigate the wavelet behaviour, we propose to separate the Real and Imaginary parts of $s^{(deO)}(t)$, introducing new functions *rpc*(.) and *ipc*(.)

$$\Re\psi^{(deO)}(t) = \Re\left\{s^{(deO)}(t-\frac{1}{2})\right\}, \Im\{\psi^{(deO)}(t)\} = \Im\left\{s^{(deO)}(t-\frac{1}{2})\right\}. \quad (18)$$

**Proposition 2**. Let $\Delta w^{(+1)} := w_0 + B$ and $\Delta w^{(-1)} := w_0 - B$; $\Delta\theta^{(+1)} := Bt_0 + w_0 t_0 + \theta_0$ and $\Delta\theta^{(-1)} := Bt_0 - w_0 t_0 - \theta_0$ be auxiliary parameters. Then

$$rpc(t) = \frac{1}{2\pi}\frac{-t_0 \cdot \sum_{i\in\{-1,+1\}} sen\,\Delta\theta^{(i)}\cos\Delta w^{(i)} t + t \cdot \sum_{i\in\{-1,+1\}}(i)\cos\Delta\theta^{(i)}\,sen\,\Delta w^{(i)} t}{t^2 - t_0^2},$$

$$ipc(t) = \frac{1}{2\pi}\frac{-t_0 \cdot \sum_{i\in\{-1,+1\}} sen\,\Delta\theta^{(i)}\,sen\,\Delta w^{(i)} t + t \cdot \sum_{i\in\{-1,+1\}}(-i)\cos\Delta\theta^{(i)}\cos\Delta w^{(i)} t}{t^2 - t_0^2}. \quad \square$$

**Proof**. Follows from trigonometry identities.

At this point, an alternative notation $rpc(t) = rpc\left(t;\Delta w^{(+1)},\Delta w^{(-1)},\Delta\theta^{(+1)},\Delta\theta^{(-1)}\right)$ and $ipc(t) = ipc\left(t;\Delta w^{(+1)},\Delta w^{(-1)},\Delta\theta^{(+1)},\Delta\theta^{(-1)}\right)$ can be introduced to explicit the dependence on these new parameters. Handling apart the real and imaginary parts of $s^{(deO)}(t)$, we arrive at

$$\sqrt{2\pi}\Re\left(s^{(deO)}(t)\right) = rpc\left(t;\pi(1+\alpha),\pi(1-\alpha),0,\frac{\pi}{2}\right) + \quad (19)$$

$$rpc(t;2\pi(1-\alpha),\pi(1+\alpha),0,0) + rpc\left(t;2\pi(1+\alpha),2\pi(1-\alpha),\frac{\pi}{2},0\right).$$

Applying now Proposition 2, after many algebraic manipulations:

$$\sqrt{2\pi}\Re\left(s^{(deO)}(t)\right) =$$

$$= \frac{1}{2}\{H_{2(1-\alpha)}(t) - H_{(1+\alpha)}(t) + M_{1+\alpha}^{1-\alpha}(t) + M_{2(1-\alpha)}^{2(1+\alpha)}(t)\}, \quad (20)$$

and $\Re\left(\psi^{(deO)}(t)\right) = \Re\left(s^{(deO)}(t-1/2)\right)$.

The analysis of the imaginary part can be done in a similar way.

**Definition 3**. (Special functions); $\nu$ is a real number, $\overline{H}_\nu(t) := \nu\frac{\cos(\nu\pi t)}{\nu\pi t}$, $0 \leq \nu \leq 1$, and

$\overline{M}_{\nu_2}^{\nu_1}(t) := \frac{1}{\pi}\frac{2|\nu_1 - \nu_2|}{1-[2t(\nu_1-\nu_2)]^2}\{sin\,\pi\nu_1 t - 2(\nu_1-\nu_2)t.\cos\pi\nu_2 t\}$ $\square$

The imaginary part of the wavelet can be found *mutatus mutandi*:

$$\sqrt{2\pi}\Im\left(s^{(deO)}(t)\right) =$$

$$= \frac{1}{2}\{\overline{H}_{2(1-\alpha)}(t) - \overline{H}_{(1+\alpha)}(t) + \overline{M}_{1+\alpha}^{1-\alpha}(t) + \overline{M}_{2(1-\alpha)}^{2(1+\alpha)}(t)\}, \quad (21)$$

and $\Im\left(\psi^{(deO)}(t)\right) = \Im\left(s^{(deO)}(t-1/2)\right)$.

The real part (as well as the imaginary part) of the complex wavelet $\psi^{(deO)}(t)$ are plotted in Figure 9, for $\alpha$ =0.1, 0.2 and 1/3.

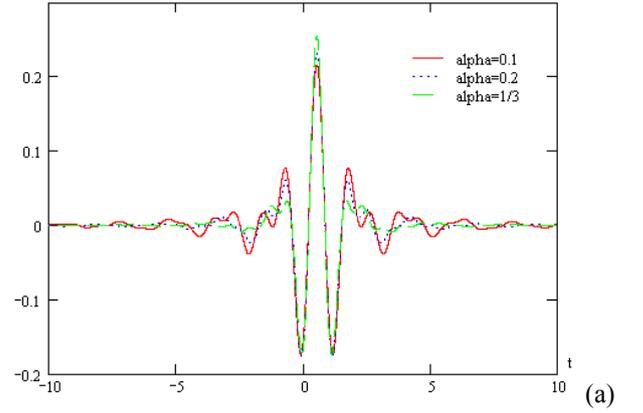

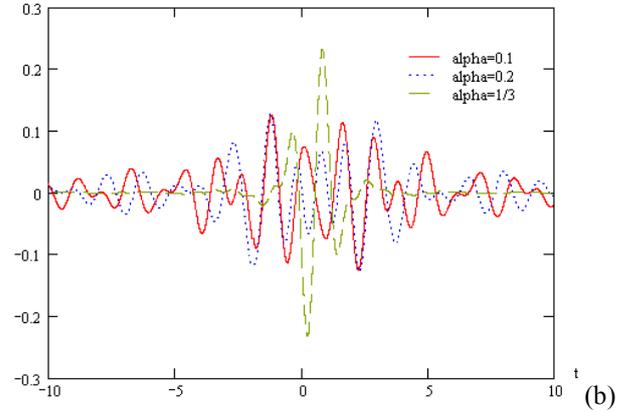

**Figure 9**. Wavelet $\psi^{(deO)}(t)$: (a) real part of the wavelet and (b) imaginary part of the wavelet. (Sketches for $\alpha$ = 0.1, 0.2 and 1/3). The effective support of such wavelets is the interval [-12,12].

## 4. Further Generalisations

Generally, the approach presented in the last section is not restricted to raised cosine filters.

**Algorithm of MRA Construction**. Let $P(w)$ be a real band-limited function, $P(w)=0$ $w>2\pi$, which satisfies the vestigial side band symmetry condition, i.e.,

$$\{P(w) + P(w-2\pi)\} = \frac{1}{2\pi} \text{ for } |w| < \pi, \quad (22)$$

then the scaling function $\Phi(w) = \sqrt{P(w)}$ defines an orthogonal MRA. $\square$

**Proposition 3**. If $P(w;\alpha)$ is a Nyquist filter of roll-off $\alpha$, and $\lambda(\alpha)$ is an arbitrary probability density function, $0<\alpha<1$, then the scaling function $\Phi(w) = \sqrt{\int_0^1 \lambda(\alpha)P(w;\alpha)d\alpha}$ defines an orthogonal MRA.

**Proof**. It is enough to show that $\sum_{n\in Z}|\Phi(w+2\pi n)|^2 = \frac{1}{2\pi}$. Given an integer $n$, then $\Phi(w+2\pi n) = \sqrt{\int_0^1 \lambda(\alpha)P(w+2\pi n;\alpha)d\alpha}$. Taking the square of both members and adding equations for each integer $n$,

$$\sum_{n \in Z} |\Phi(w+2\pi n)|^2 = \int_0^1 \lambda(\alpha) \sum_{n \in Z} P(w+2\pi n;\alpha) d\alpha$$ . Since $P(w;\alpha)$ is a Nyquist filter, $\sum_{n \in Z} P(w+2\pi n;\alpha) = \frac{1}{2\pi}$ and the proof follows. □

The most interesting case of such generalisation corresponds to a "weighting" of square-root-raise-cosine filter.

**Corollary**. (Generalised raised-cosine MRA). If $P(w;\alpha)$ is the raised cosine spectrum with a (continuous) roll-off parameter $\alpha$, i.e.,

$$P(w;\alpha) := \begin{cases} \frac{1}{2\pi} & 0 \leq |w| < (1-\alpha)\pi \\ \frac{1}{4\pi}\left\{1+\cos\frac{1}{2\alpha}(|w|-\pi(1-\alpha))\right\} & (1-\alpha)\pi \leq |w| < (1+\alpha)\pi \\ 0 & |w| \geq (1+\alpha)\pi \end{cases}$$

and $\lambda(\alpha)$ is an arbitrary probability density function defined over the interval $0<\alpha<1$, the scaling function $\Phi(w) = \sqrt{\int_0^1 \lambda(\alpha)P(w;\alpha)d\alpha}$ defines an orthogonal MRA. □

## 5. Implementation on MATLAB<sup>TM</sup> and a few Applications

Aiming to investigate some potential applications of such wavelets, software to compute them should be written. Nowadays one of the most powerful software supporting wavelet analysis is the Matlab[TM] [24], especially when the wavelet graphic interface is available. In the Matlab[TM] wavelet toolbox, there exist five kinds of wavelets (type the command *waveinfo* on the prompt): (*i*) crude wavelets (*ii*) Infinitely regular wavelets (*iii*) Orthogonal and compactly supported wavelets (*iv*) biorthogonal and compactly supported wavelet pairs. (*v*) complex wavelets. Complex wavelets present as a rule interesting symmetry and possess closed expression. In contrast with standard continuous complex wavelets (cgau*N* Gaussian, cmor*F*b-*F*c Morlet, Shannon shan*F*b-*F*c, fbsp*M*-*F*b-*F*c frequency *B*-spline), the "de Oliveira" scaling function do exist and these wavelets are orthogonal so the reconstruction property is insured.

The m-files to allow the computation of de Oliveira wavelet transform are currently (*freeware*) available at the URL: http://www2.ee.ufpe.br/codec/WEBLET.html (new wavelets). The family de Oliveira includes: deo (short name) and cdeo (short name). Support width infinite, effective support [-12 12]. Since $0 \leq \alpha \leq 1/3$, the following five default values were assumed: $\alpha = 0, 1/15, 2/15, 1/4, 4/15$ and $1/3$. For instance, the complex wavelet associated with $\alpha = 1/3$ is denoted cdeo0.33333. A naive application is presented in the sequel so as to illustrate the potential of this tool.

The continuous de Oliveira wavelet transform ($a \neq 0$) of an arbitrary signal $f(t)$ can be computed according to:

$$C_{a,b} := \frac{1}{\sqrt{|a|}} \int_{-\infty}^{+\infty} f(t) \psi^{(deO)}\left(\frac{t-b}{a}\right) dt, \qquad (23)$$

where $\psi^{(deO)}(.)$ is defined by Eqn(18), (20), (21).

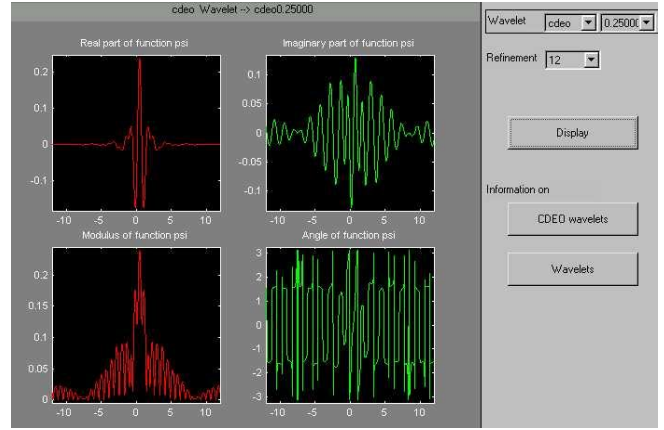

**Figure 10**. Complex orthogonal "de Oliveira" wavelet display over Matlab[TM] using the *wavemenu* command.

### 5.1 Fault Detection in Transmission Lines

Just to exemplify the potential of this tool in signal processing, a few signal derived by fault simulation in a 138kV transmission line by the ATP (Alternative Transient Program) were analysed [25]. Figure 11 shows a naive scheme for the section of the transmission system used in simulations. Here, $Z_{LT}$ is the transmission line impedance between A / B; $Z_{e1}$ and $Z_{e2}$ are the Thevenin equivalent impedance in the terminal A and B, respectively; $E_{e1}$ and $E_{e2}$ denote the voltage source de of the Thevenin equivalent circuit in the two terminals. All faults are derived at a sampling rate of 128 samples/cycle.

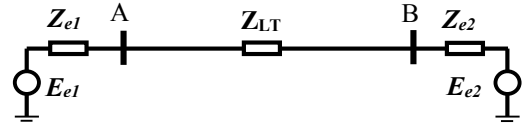

**Figure 11**: Simplified diagram of the three-phase line adopted in ATP simulations.

The faults were simulated by ATP, considering that the power quality monitoring took place at terminal A. The transmission line is assumed to be totally transposed, with a distributed parameter model. The total simulation slot was eight cycles, and all faults (with zero fault impedance) were simulated four cycles after the start of simulation.

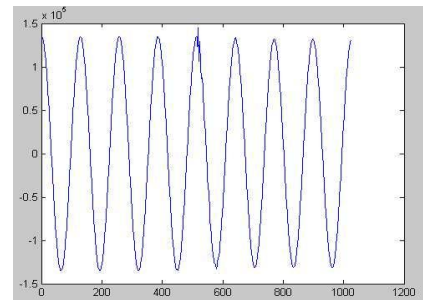

**Figure 12**. Voltage sample of a single-phase after a phase-phase fault simulated by ATP for the transmission line shown in Figure 11. A short transient occurred at four cycles after starting simulation.

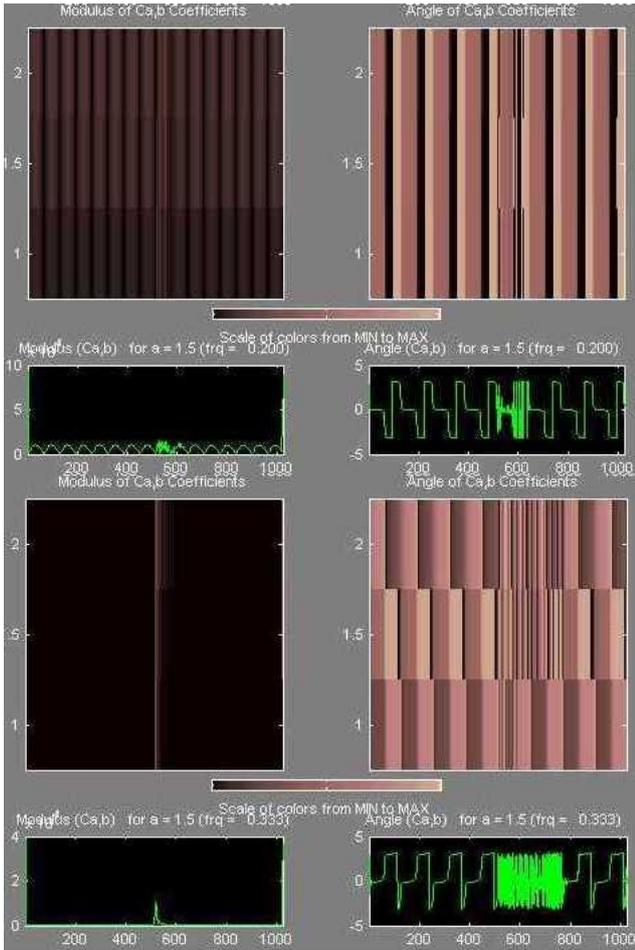

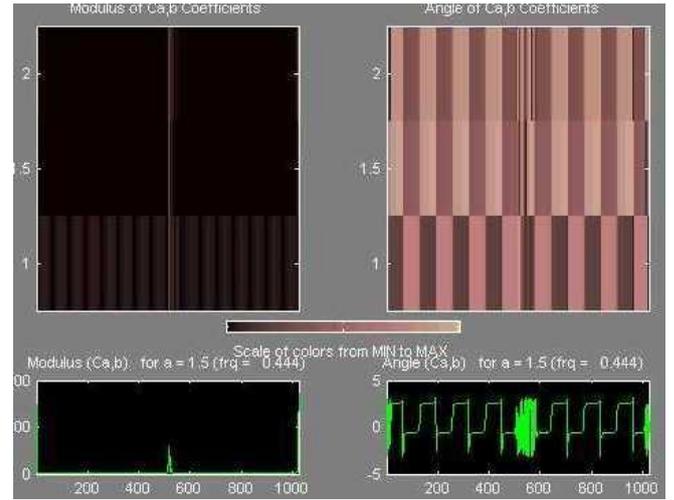

**Figure 14**. Modulus of $C_{ab}$ and coefficient line $a$=1 (left), Phase of $C_{ab}$ and coefficient line $a$=1 (right) for the signal (Fig. 12). Scale settings: Min.1, step 0.5, Max. 2; colour-map: "pink", number of colour: 256. Different complex analysing wavelets: (a) cdeo0.08333 (b) cdeo0.33333 .

## 6. Conclusions

This paper introduced a new family of complex orthogonal wavelets, which was derived from the classical Nyquist criterion for ISI elimination in Digital Communication Systems. Properties of both the scaling function and the mother wavelet were investigated. This new wavelet family can be used to perform an orthogonal Multiresolution Analysis. A new function termed *PCOS* was introduced, which is offered as a powerful tool in matters that concern raised cosines. An algorithm for the construction of MRA based on vestigial side band filters was presented. A generalisation of the (square root) raised cosine wavelet was also proposed yielding a broad class of orthogonal wavelets and MRA. These wavelets have been implemented on Matlab™ (wavelet toolbox) and an application to fault detection in transmission lines described. These new wavelets seem to be particularly suitable as natural candidates to replace *Sinc*(.) pulses in standard OFDM systems. Our group is currently investigating this topic.

## ACKNOWLEDGMENTS

The authors thank Dr. Ricardo Campello de Souza for his countless constructive criticism. The first author is glad in contributing to this memorial tribute to Dr. Gerken, who he had the privilege to share a pleasant companionship.

**Figure 13**. Modulus of $C_{ab}$ and coefficient line $a$=1 (left), Phase of $C_{ab}$ and coefficient line $a$=1 (right) for the signal (Fig. 12). Scale settings: Min.1, step 0.5, Max. 2; colour-map: "pink", number of colour: 256. Different complex analysing wavelets: (a) cgau1 (b) cgau2.

The wavelet analysis of the single-phase voltage signal shown in Fig.12 was carried out applying different complex wavelets as an attempt to locate the position where the fault start. Figure 13 and 14 illustrate the results of magnitude and phase-angle jump for two examples: "complex gaussian" and "de Oliveira" wavelets. It can be observed that cdeo can furnish nice time localisation properties.

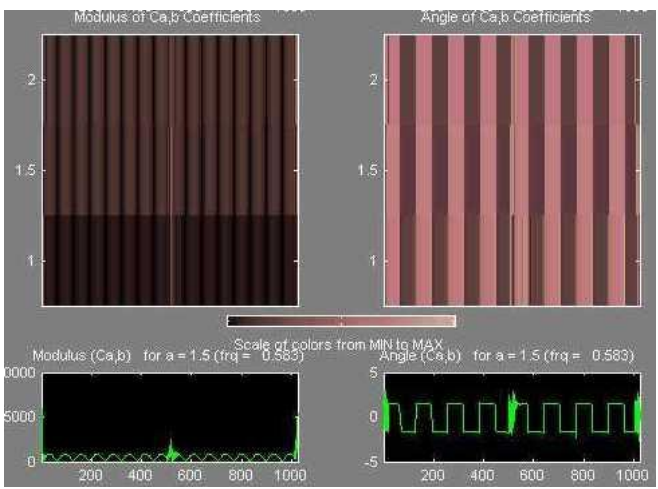

**Hélio Magalhães de Oliveira** was born in Arcoverde, PE. He received both the B.Eng. and the M.Sc. degrees in Electrical Engineering from the Federal University of Pernambuco (UFPE) Brazil, in 1980 and 1983, respectively. Then he joined the staff of the Electronics and Systems Department at the same University as a lecturer. In 1992, he earned the "Docteur de l'Ecole Supérieure des Télécommunications" degree, in Paris, France. He was elected honoured professor by twenty electrical engineering undergraduate groups and chosen as the godfather of two engineering graduation. His publications are available (in the .pdf format) at http://www2.ee.ufpe.br/codec/publicacoes.html. He was the head of the UFPE electrical engineering graduate program from 1992 to 1996. Dr. de Oliveira authored the book *Análise de Sinais para Engenheiros: Uma Abordagem via Wavelets*, São Paulo: Editora Manole, 2003 (in press). His current research interests include: digital signal processing, wavelets, data communication, applied information theory with emphasis on error-control coding. Dr. de Oliveira is a member of the Brazilian Telecommunication Society (SBrT) and the Institute of Electrical and Electronics Engineers (IEEE).

**Luciana Reginaldo Soares** was born in Rio de Janeiro-RJ, Brazil, in 1973. She received both the B.Eng. and the M.Sc. degrees in Electrical Engineering from the Federal University of Pernambuco (UFPE), Brazil, in 1997 and 2001, respectively. She is currently a doctoral student at the same University. Her research interests include digital signal processing, wavelets and power system analysis. L.R. Soares is a student member of the "Conseil International des Grands Réseaux Electriques" (CIGRÉ).

**Tiago Henrique Falk** was born in Recife, Pernambuco, in September 1979. He was recipient of CNPq undergraduate scholarship and he was awarded with the "Prof. Newton Maia Young Scientist Prize" (UFPE-CNPq) in 2001. He received the B.Eng. degree in Electrical engineering from the Federal University of Pernambuco, in 2002. He enrolled in Queen's University, Canada, where he is presently working toward a Ph.D. His interest areas are: Communication Theory and Signal Processing. Tiago Falk is student member of both the Institute of Electrical and Electronic Engineering (IEEE) and the Brazilian Telecommunication Society (SBrT).